\magnification\magstep1
\baselineskip = 18pt

\def\square{\vcenter{\hrule height1pt
\hbox{\vrule width1pt height4pt \kern4pt
\vrule width1pt}
\hrule height1pt}}
\def \Bbb {\bf}
\centerline{Convex Bodies with Few Faces}\medskip
\centerline{by}\medskip
\centerline{Keith Ball$^{(1)}$}
\centerline{Texas A\&M University}
\centerline{College Station, TX  77843}\bigskip
\centerline{and}\bigskip
\centerline{Alain Pajor}
\centerline{U.E.R. de Math\'ematiques}
\centerline{Universit\'e de Paris VII}
\centerline{2 Place Jussieu}
\centerline{75251 PARIS  CEDEX 05}\bigskip

\noindent {\bf Abstract.} It is proved that if $u_1,\ldots, u_n$ are
vectors in ${\Bbb R}^k, k\le n, 1 \le p < \infty$ and

$$r = \bigg( {1\over k} \sum ^n_1 |u_i|^p\bigg)^{1\over p}$$

\noindent then the volume of the symmetric convex body whose boundary
functionals are \hfil\break $\pm u_1,\ldots, \pm u_n$, is bounded from
below as

$$|\{ x\in {\Bbb R}^k\colon \ |\langle x,u_i\rangle | \le 1 \ \hbox{for
every} \ i\}|^{1\over k} \ge {1\over \sqrt{\rho}r}.$$

\noindent An application to number theory is stated.\vskip1.5in

\noindent A.M.S. (1980) Subject Classification: \ 52A20, 10E05

\noindent $^{(1)}$Partially supported by N.S.F. DMS-8807243.
\vfill\eject

\noindent {\bf \S 0. Introduction.}

In [V], Vaaler proved that if $Q_n = [-{1\over 2}, {1\over 2}]^n$ is the
central unit cube in ${\Bbb R}^n$ and $U$ is a subspace of ${\Bbb R}^n$
then the volume $|U \cap Q_n|$, of the section of $Q_n$ by $U$ is at least
1. This result may be reformulated as follows: \ if $u_1,\ldots, u_n$ are
vectors in ${\Bbb R}^k, 1 \le k \le n$ whose Euclidean lengths satisfy
$\sum\limits ^n_1 |u_i|^2 \le k$ then

$$|\{x\in {\Bbb R}^k\colon \ |\langle x,u_i\rangle| \le 1 \ {\rm for \
every} \ i\}|^{1\over k} \ge 2.$$

A related theorem, (Theorem 1, below) in which the condition $\sum |u_i|^2
\le k$ is replaced by $\max\limits _i|u_i| \le 1$ was proved by Carl and
Pajor [C-P] and Gluskin [G]. Gluskin's methods enable him to obtain sharp
results in limiting cases which in turn have applications in harmonic
analysis. Results closely related to Theorem 1 were also obtained by
B\'ar\'any and F\"uredi [B-F] and Bourgain, Lindenstrauss and Milman
[B-L-M].\medskip

\noindent {\bf Theorem 1.} There is a constant $\delta >0$ so that if
$u_1,\ldots, u_n \in {\Bbb R}^k, 1 \le k \le n$ are vectors of length at
most 1 then

$$|\{x\in {\Bbb R}^k\colon \ |\langle x,u_i\rangle | \le 1 \ \hbox{for
every}\ i\}|^{1\over  k} \ge {\delta \over \sqrt{1+\log {n\over k}}}.$$

The estimate is best possible if $n$ is at most exponential in $k$, apart
from the value of the constant $\delta$. This is demonstrated by an example
which had appeared some time earlier in a paper of Figiel and Johnson,
[F-J]. Theorem 1 gives a lower bound on the volume ratios of the unit
balls
of $k$-dimensional subspaces of $\ell^n_\infty$ and hence on the distance
of these subspaces from Euclidean space.

Regarding Theorem 1 as a ``$p=\infty$'' version of Vaaler's ``$p=2$''
result, Kashin asked whether a similar result holds for $2<p<\infty$.
This question is answered in the affirmative by the following theorem.
\medskip

\noindent {\bf Theorem 2.} Suppose $u_1,\ldots, u_n \in {\Bbb R}^k$ with
$k\le n, 1 \le p <\infty$ and let $r = ({1\over k} \sum\limits ^n_1
|u_i|^p)^{1\over p}$. Then

$$|\{x \in {\Bbb R}^k\colon \ |\langle x,u_i)| \le 1 \ {\rm for \
every} \ i\}|^{1\over k}
 \ge \left\{\matrix{{2\sqrt{2}\over \sqrt{p}r}&{\rm if}& p\ge
2\hfill\cr {1\over r}&{\rm if}&1\le p \le 2.\hfill\cr}\right.$$

\noindent The lower bound is best possible (up to a constant) provided
$e^pk \le n\le e^k$.

\noindent {\bf Remark.} The slightly stronger result for $p\ge 2$ is
isolated since for $p=2$ it gives back exactly Vaaler's result.
\medskip

Theorem 1 follows immediately from Theorem 2 by a standard optimisation
argument. If $(u_i)^n_1$ in ${\Bbb R}^k$ all have norm at most 1 then for
any $p\in [1,\infty)$,

$$\bigg({1\over k} \sum ^n_1 |u_i|^p\bigg)^{1\over p} \le \Big({n\over
k}\Big)^{1\over p}$$

\noindent so that

$$|\{x\colon \ |\langle x,u_i\rangle| \le 1 \ \hbox{for every} \
i\}|^{1\over k} \ge {2\sqrt{2}\over \sqrt{p} ({n\over k})^{1\over p}}$$

\noindent (for $p\ge 2$) and the latter is at least ${2\over \sqrt{e}
\sqrt{1+\log {n\over k}}}$ when $p = 2(1+\log {n\over k})$.

With the careful use of well-known methods for estimating the entropy of
convex bodies it is possible to obtain more general (but less precise)
estimates than that provided by Theorem 2; (see [B-P]). The purpose of this
paper is to provide a very short proof of Theorem 2 and, a fortiori,
Theorem 1.

Vaaler originally proved his theorem because of its applications to the
geometry of numbers. The last section of this paper includes a statement of
the generalisation of Siegel's lemma which follows from Theorem 2.
\vfill\eject

\noindent {\bf \S 1. The lower bound.}

The proof of Theorem 2 makes use of the following result from [Me-P] which
was designed to extend Vaaler's theorem in a different direction: \ it
estimates the volumes of sections of the unit balls of the spaces
$\ell^n_p, 1 \le p \le \infty$. For $1 \le p \le \infty, n \in {\Bbb N}$
let

$$B^n_p = \bigg\{x\in {\Bbb R}^n\colon \ \sum^n_1 |x_i|^p \le 1\bigg\}$$

\noindent be the unit ball of $\ell^n_p$.\medskip

\noindent {\bf Theorem 3.} Let $U$ be a $k$-dimensional subspace of ${\Bbb
R}^n$; if $1 \le p \le q \le \infty$ then

$${|B^n_p \cap U|\over |B^k_p|} \le {|B^n_q \cap U|\over |B^k_q|}.\eqno
\square$$

\noindent {\bf Remark.} The case $p =2, q=\infty$ is Vaaler's theorem since
then, the left side is 1 and the inequality states that

$$|B^n_\infty \cap U| \ge |B^k_\infty| \ge 2^k.$$

For notational convenience, the proof of Theorem 2 is divided into several
short lemmas. The first is no more than a convenient form of H\"older's
inequality. For $k \in {\Bbb N}, S^{k-1}$ will denote the Euclidean sphere
in ${\Bbb R}^k$ and $\sigma = \sigma_{k-1}$, the rotationally invariant
probability measure on $S^{k-1}$. Also let $v_k$ be the volume of the
Euclidean unit ball in ${\Bbb R}^k$.\medskip

\noindent {\bf Lemma 4.} Let $C$ and $B$ be symmetric convex bodies in
${\Bbb R}^k$ with Minkowski gauges $\|\cdot\|_C$ and $\|\cdot \|_B$
respectively. Then for $p>0$

$$\Big({|C|\over |B|}\Big)^{1\over k} \ge \bigg({k+p\over k|B|} \int_B
\|x\|^p_C dx\bigg)^{-{1\over p}}.$$

\noindent {\bf Proof.}

$$\eqalignno{\Big({|C|\over |B|}\Big)^{1\over k} &= \bigg( {v_k\over |B|}
\int_{S^{k-1}} \|\theta\|^{-k}_C d\sigma(\theta)\bigg)^{1\over k}\cr
&= \bigg({kv_k\over |B|} \int_{S^{k-1}} \Big({\|\theta\|_B\over
\|\theta\|_C}\Big)^k \int ^{\|\theta\|^{-1}_B}_0 r^{k-1} dr
d\sigma(\theta)\bigg)^{1\over k}\cr
&= \bigg( {1\over |B|} \int_B \Big({\|x\|_B\over \|x\|_C}\Big)^k
dx\bigg)^{1\over k}\cr
&\ge \bigg({1\over |B|} \int_B \Big({\|x\|_B\over \|x\|_C}\Big)^{-p}
dx\bigg)^{-{1\over p}}\cr
&= \bigg( {k+p\over k|B|} \int_B \|x\|^p_C dx\bigg)^{-{1\over
p}}.&\square}$$

\noindent (Lemma 4 appears in [Mi-P] as Corollary 2.2.)\medskip

\noindent {\bf Lemma 5.} Suppose $u_1,\ldots, u_n \in {\Bbb R}^k$ with
$k\le n$ and $1\le p < \infty$. Then

$$\eqalign{&|\{x\in {\Bbb R}^k\colon \ |\langle x_i,u_i\rangle| \le 1 \
\hbox{for every} \ i\}|^{1\over k}\cr
&\qquad \ge 2\bigg({k+p\over k} \sum ^n_{i=1} {1\over |B^k_p|} \int_{B^k_p}
|\langle x,u_i\rangle|^p dx\bigg)^{-{1\over p}}.}$$

\noindent {\bf Proof.}	Define $T\colon \ {\Bbb R}^k \to {\Bbb R}^n$ by
$(Tx)_i = \langle x,u_i\rangle, 1 \le i \le n$ and let $U = T({\Bbb R}^k)$.
The problem is to estimate from below

$$|T^{-1}(B^n_\infty)|^{1\over k} = |T^{-1}(U\cap B^n_\infty)|^{1\over
k}.$$

\noindent By Theorem 3,

$$\eqalignno{|U\cap B^n_\infty|^{1\over k} &\ge 2 \Big({|U\cap B^n_p|\over
|B^k_p|}\Big)^{1\over  k}\cr
\noalign{\hbox{and so}}
|T^{-1}(B^n_\infty)|^{1\over k} &\ge 2\Big({|T^{-1}(B^n_p)|\over
|B^k_p|}\Big)^{1\over k}.}$$

\noindent Regard $T$ as an operator: \ $\ell^k_p \to \ell^n_p$. Then by
Lemma 4,

$$\eqalignno{2\Big({|T^{-1}(B^n_p)|\over |B^k_p|}\Big)^{1\over k} &\ge 2
\bigg( {k+p\over k|B^k_p|} \int _{B^k_p} \|Tx\|^p dx\bigg)^{-{1\over p}}\cr
&= 2 \bigg({k+p\over k|B^k_p|} \int_{B^k_p} \sum^n_1 |\langle
x,u_i\rangle|^pdx\bigg)^{-{1\over p}}.&\square}$$

\noindent {\bf Proof of Theorem 2.} Let $(u_i)^n_1$ and $p$ be as above.
For each $i$ let $v_i$ be the unit vector in the direction of $u_i$. By
Lemma 5,

$$\eqalign{&|\{x\in {\Bbb R}^k\colon \ |\langle x,u_i\rangle| \le 1 \
\hbox{for every} \ i\}|^{1\over k}\cr
&\ge 2 \bigg( {k+p\over k} \sum ^n_{i=1} {1\over |B^k_p|} \int _{B^k_p}
|\langle x,u_i\rangle| ^p dx\bigg)^{-{1\over p}}\cr
&= 2 \bigg( {k+p\over k} \sum ^n_1 |u_i|^p \cdot {1\over |B^k_p|}
\int_{B^k_p} |\langle x,v_i\rangle|^p dx\bigg)^{-{1\over p}}\cr
&\ge 2 \bigg( {1\over k} \sum^n_{i=1} |u_i|^p\bigg)^{-{1\over p}} \min
\bigg({k+p\over |B^k_p|} \int_{B^k_p} |\langle
x,v\rangle|^pdx\bigg)^{-{1\over p}}}$$

\noindent where the minimum is taken over all vectors $v$ of Euclidean
length 1. So to complete the proof it suffices to show that for such a
vector $v$,

$$\bigg( {k+p\over |B^k_p|} \int_{B^k_p} |\langle x,v\rangle|^p
dx\bigg)^{1\over p} \le \left\{\matrix{\sqrt{p\over 2}&{\rm if}&p\ge
2\hfill\cr 2&{\rm if}&1\le p<2.\hfill\cr}\right.$$

\noindent Let $(x^{(j)})^k_1$ and $(v^{(j)})^k_1$ be the coordinates of the
vectors $x$ and $v$ in ${\Bbb R}^k$. For $p\ge 2$, observe that the
functions $(x^{(j)} v^{(j)})$ on $B^k_p$ form a conditionally symmetric
sequence, so by Khintchine's inequality and H\"older's inequality (for
$\sum\limits ^n_1 v^{(j)2} = 1$),

$$\eqalign{&\bigg({k+p\over |B^k_p|} \int_{B^k_p} \bigg|\sum ^k_1 x^{(j)}
v^{(j)}\bigg|^p dx\bigg)^{1\over p}\cr
&\qquad \le \sqrt{p\over 2} \bigg({k+p\over |B^k_p|} \int_{B^k_p}
\bigg(\sum ^k_1 x^{(j)2} v^{(j)2} \bigg)^{p\over 2} dx\bigg)^{1\over p}\cr
&\qquad \le \sqrt{p\over 2} \bigg( {k+p\over |B^k_p|} \int_{B^k_p} \sum
^n_1 |x^{(j)}|^p \cdot v^{(j)2} dx\bigg)^{1\over p}\cr
&\qquad = \sqrt{p\over 2} \bigg( {k+p\over |B^k_p|} \int_{B^k_p}
|x^{(1)}|^p dx\bigg)^{1\over p}\cr
&\qquad = \sqrt{p\over 2}.}$$

\noindent For $1\le p < 2$ it is easily checked that

$$\eqalign{&\bigg({k+p\over |B^k_p|} \int_{B^k_p} |\langle x,v\rangle|^p
dx\bigg)^{1\over p}\cr
&\qquad \le (k+p)^{1\over p} \bigg({1\over |B^k_p|} \int_{B^k_p} \langle
x,v\rangle^2 dx\bigg)^{1\over 2}\cr
&\qquad = (k+p)^{1\over p} \bigg({1\over |B^k_p|} \int_{B^k_p} (x^{(1)})^2
dx\bigg)^{1\over 2}}$$

\noindent and the last expression can be (rather roughly) estimated by 2
using standard inequalities involving logarithmically concave
functions.$\hfill \square$\medskip

\noindent {\bf Remark.} The proof of Theorem 2 can be simplified even
further if the integration over $B^k_p$ is replaced by integration over
$S^{k-1}$ (and H\"older's inequality applied here). The proof was presented
as above because Lemma 5 has some intrinsic interest: \ for example it may
be used to recover Gluskin's precise estimate as follows. Suppose $m\in
{\Bbb N}$ and the vectors $(z_i)^m_1 \in {\Bbb R}^k$ satisfy

$$|z_i| \le \Big(\log \Big( 1 + {m\over k}\Big)\Big)^{-{1\over 2}}, \quad 1
\le i \le m.$$

\noindent For $\varepsilon >0$, let $W(\varepsilon)$ be the set

$$\Big\{ x\in {\Bbb R}^k\colon	\ \max_j|x^{(j)}|\le 1, \max_i |\langle
x,z_i\rangle| \le {1\over \varepsilon}\Big\}\colon$$

\noindent that is, $W(\varepsilon)$ is the intersection of the cube
$B^k_\infty$ with $m$ ``bands'' of width at most \break ${2\over
\varepsilon}\sqrt{\log(1+{m\over k})}$.  Then $|W(\varepsilon)|^{1\over
k}
\to 2$ as $\varepsilon \to 0$, uniformly in $k$ and $m$. To see this, apply
Lemma 5 with $n= k+m$, the first $k, u_i$'s being the standard basis
vectors of ${\Bbb R}^k$ and the remaining $m$ being the vectors
$(\varepsilon z_i)^m_1$. If $e_j$ is a standard basis vector,

$${k+p\over |B^k_p|} \int_{B^k_p} |\langle x,e_j\rangle|^pdx = 1$$

\noindent and so Lemma 5 (and the proof of Theorem 2) show that for each
$p\ge 2, |W(\varepsilon)|^{1\over k} \ge \break 2(1 + {m\over k}({p\over
2})^{p\over 2} \varepsilon^p(\log(1 + {m\over k}))^{-{p\over 2}})^{-{1\over
p}}$ and the latter is at least ${2\over 1+\sqrt{e}\cdot \varepsilon}$ if
$p = \max(2, \break 2\log(1 + {m\over k}))$.$\hfill\square$\medskip

As was briefly mentioned earlier, more general estimates than that of
Theorem 2 are obtained in [B-P] (for entropy numbers instead of volumes).
It is worth noting however that even the argument of Theorem 2 can be used
to give the following: \ there is a constant $c$ so that if $u_1,\ldots,
u_n \in {\Bbb R}^k, k\le n$ and $T\colon \ \ell^k_2 \to \ell^n_\infty$ is
given by $(Tx)_i = \langle u_i,x\rangle, \ 1 \le i \le n$, then the $k^{\rm
th}$ entropy number of $T$ satisfies

$$e_k(T) \le {c\sqrt{p}\over \sqrt{k}} \bigg({1\over k} \sum^n_1
|u_i|^p\bigg)^{1\over p} \Big(1+\log {n\over k}\Big)^{1\over p}$$

\noindent and hence

$$e_k(T) \le {ec\over \sqrt{k}} \sqrt{1+\log {n\over k}}\cdot \|T\|$$

\noindent (taking $p = 2(1+\log {n\over k}))$. To obtain this one uses
Sch\"utt's estimates, [5], for the entropy numbers of the formal identity
from $\ell^n_p$ to $\ell^n_\infty$ in place of the result of Meyer and
Pajor, and the dual Sudakov inequality of Pajor and Tomczak, [P-T] in place
of the application of H\"older's inequality.\vfill\eject

\noindent {\bf \S 2. An application to linear forms.}

As stated in the introduction, Vaaler's original result has applications to
the geometry of numbers. One such, a sharpened form of Siegel's lemma, is
given in [B-V]. Using the arguments of Bombieri and Vaaler and Theorem 2,
one can obtain the generalisation of their result, contained in Theorem 6,
below. Some notation is needed. If $A$ is  a $k\times n$ matrix of reals
with independent rows $(1\le k \le n)$, denote by $v_j	= v_j(A), \ 1\le j
\le k$, the rows of $A$. Let $(e_i)^n_1$ be the standard basis of ${\Bbb
R}^n$ and denote by $c_i$, the distance (in the Euclidean norm) of $e_i$
from the span of the $v_j$'s in ${\Bbb R}^n$. (So if $A_i$ is the matrix
with $k+1$ rows, $v_1,\ldots, v_n,e_i$ then

$$c^2_i = {\det (A_iA^*_i)\over \det (AA^*)} \quad {\rm for}\quad 1 \le i
\le n.)$$

\noindent {\bf Theorem 6.} Let $A$ be a $k\times n$ matrix with rank $k$
and integral entries. With the notation above, the system $Ax=0$ admits
$n-k$ linearly independent solutions

$$z^{(r)} = (z^{(r)}_1,\ldots, z^{(r)}_n) \in {\Bbb Z}^n,\ 1 \le r \le
n-k$$

\noindent so that for every $p\ge 2$,

$$\prod_{1\le r \le n-k} \max_i |z^{(r)}_i| \le D^{-1} \sqrt{p\over 2}
\bigg( {1\over n-k} \sum^n_1 c_i^p\bigg)^{1\over p} \sqrt{\det AA^*}$$

\noindent where $D$ denotes the G-C-D of all $k\times k$ determinants
extracted from $A$.$\hfill \square$

\noindent {\bf Remark.} The principal importance of such a generalisation
of Bombieri and Vaaler's result is that it takes into account, more
strongly, the form of the matrix $A$. If the $c_i$'s are all about the same
size, then for $p>2$, the expression

$$\bigg({1\over n-k} \sum c_i^p\bigg)^{1\over p}$$

\noindent is small compared with the corresponding expression in which $p$
is replaced by 2.\vfill\eject

\noindent {\bf References.}

\item{[B-F]} I. B\'ar\'any and Z. F\"uredi, Computing the volume is
difficult, Discrete Comput. Geom. 2 (1987), 319-326.

\item{[B-L-M]} J. Bourgain, J. Lindenstrauss and V.D. Milman, Approximation
of zonoids by zonotopes, Acta Math. 162 (1989), 73-141.

\item{[B-P]} K.M. Ball and A. Pajor, On the entropy of convex bodies with
``few'' extreme points, in preparation.

\item{[B-V]} E. Bombieri and J. Vaaler, On Siegel's lemma, Invent. Math. 73
(1983), 11-32.

\item{[C-P]} B. Carl and A. Pajor, Gelfand numbers of operators with values
in a Hilbert space, Invent. Math. 94 (1988), 479-504.

\item{[F-J]} T. Figiel and W.B. Johnson, Large subspaces of $\ell^n_\infty$
and estimates of the Gordon-Lewis constants, Israel J. Math. 37 (1980),
92-112.

\item{[G]} E.D. Gluskin, Extremal properties of rectangular parallelipipeds
and their applications to the geometry of Banach spaces, Math. Sbornik, 136
(178) (1988), 85-95.

\item{[Me-P]} M. Meyer and A. Pajor, Sections of the unit ball of
$\ell^n_p$, J. Funct. Anal. 80 (1988), 109-123.

\item{[Mi-P]} V.D. Milman and A. Pajor, Isotropic position and inertia
ellipsoids and zonoids of the unit ball of a normed $n$-dimensional space.
Israel seminar on G.A.F.A. (1987-88), Springer-Verlag, Lecture notes in
Math. \#1376 (1989).

\item{[P-T]} A. Pajor and N. Tomczak-Jaegermann, Subspaces of small
codimension of finite-dimensional Banach spaces, Proc. A.M.S. 97 (1986),
637-642.

\item{[S]} C. Sch\"utt, Entropy numbers of diagonal operators between
symmetric Banach spaces, Journal Approx. Th. 40 (1984), 121-128.

\item{[V]} J.D. Vaaler, A geometric inequality with applications to linear
forms, Pacific J. Math. 83 (1979), 543-553.

\end